\documentclass[11pt]{article}
\usepackage[english]{babel}
\usepackage{amsmath}
\usepackage{amsthm}
\usepackage{amssymb}
\usepackage[dvips]{color}
\newtheorem{theorem}{Theorem}[section]
\newtheorem{definition}[theorem]{Definition}
\newtheorem{lemma}[theorem]{Lemma}

\newtheorem{remark}[theorem]{Remark}

\newcommand{\pn}{\par\noindent}

\newcommand{\Er}{\mathcal{R}}
\newcommand{\ov}{\overline}

\newcommand{\rr}{\mathbb{R}}

\newcommand{\hh}{\mathbb{H}}

\newcommand{\pp}{\partial}

\title{\bf  Non commutative functional calculus:  unbounded operators }

\author{Fabrizio Colombo\\Dipartimento di Matematica\\Politecnico di
Milano\\Via Bonardi, 9\\20133 Milano,
Italy\\fabrizio.colombo@polimi.it \and Graziano
Gentili\\Dipartimento di Matematica\\ Universit\'a di Firenze
\\Viale Morgagni, 67 A\\
Firenze, Italy\\ gentili@math.unifi.it
\\
\and Irene Sabadini\\Dipartimento di
Matematica\\Politecnico di Milano\\Via Bonardi, 9\\20133 Milano,
Italy\\irene.sabadini@polimi.it
\and Daniele C. Struppa\\Department of Mathematics \\and Computer Sciences
\\Chapman University\\Orange, CA 92866 USA,
\\struppa@chapman.edu}

\date{ }
\begin{document}
\maketitle
\begin{abstract}
In a recent work, \cite{cgss}, we developed a functional calculus for bounded operators
defined on quaternionic Banach spaces. In this paper we show how the results from \cite{cgss}
can be extended to the unbounded case, and we highlight the crucial differences between the two
cases. In particular, we deduce a new eigenvalue equation, suitable for the construction
of a functional calculus for operators whose spectrum is not necessarily real.
 \end{abstract}



\section{Introduction}

Let $V$ be a Banach space over the skew field $\mathbb{H}$ of
quaternions, and let $E$ be the Banach space of right linear
operators acting on it. In a recent paper, \cite{cgss}, we have
shown that if $T$ is a bounded operator in $E$ the standard
eigenvalue problems (i.e. the search of the values
$\lambda\in\mathbb{H}$ for which $\lambda \mathcal{I}-T$ or
$\mathcal{I}\lambda -T$ are not invertible) do not lead to a useful
functional calculus. The main reason for this difficulty consists
in the fact that the inverse of the operator $\lambda
\mathcal{I}-T$ does not correspond, unlike what happens in the
complex case, to an integral kernel of Cauchy type for
non-commuting variables. In order to overcome this difficulty, we
introduced a totally new eigenvalue problem, which we called the
$S-$eigenvalue problem, and which, in the commutative case, is
identical with the usual problem. The key observation which
allowed the study of such an eigenvalue problem, was the recent
development of slice regular functions, \cite{gs}, \cite{advances}. These
functions are likely to be the appropriate generalization of
standard holomorphic functions, and they can be represented
through a new, non-commutative, Cauchy kernel. It is by using this
kernel, that one is able to develop a functional calculus for
bounded operators in $E.$ The specific techniques utilized in
\cite{cgss} did not allow us to tackle the important
study of unbounded operators. These operators, however, are of
great importance both in mathematics and in physics, because they
can be used to write  functions of operators, such as the
exponential of a closed operator, for which the corresponding
power series expansion is not convergent and thus not suitable to
define it. In particular, in quantum mechanics, the exponential
function of an operator defines the evolution operator associated
to Schr\"odinger equation. As it is well known, Quantum Mechanics
can be formulated in the real, complex and quaternionic setting
(see \cite{adler}), and for this reason it is important to
introduce a quaternionic version of the functional calculus to
allow the study of exponentials for quaternionic operators.

Finally, we point out that the spectral theory can be extended to the case of $n$-tuples of
operators. For the complex case, the reader is referred to \cite{Taylor}, \cite{Taylor3} and, by using the notion
of monogenic functions with values in a Clifford algebra (see \cite{bds} for the case of one variable and \cite{nostro libro} for the case of several variables), to \cite{jefferies} and the references therein.

{\it Acknowledgements} The first and third authors are grateful to Chapman University for the
hospitality during the period in which this paper was written. They are also indebted to
G.N.S.A.G.A. of INdAM and
the Politecnico di Milano for partially supporting their visit.

\section{ Preliminary results}

\subsection{Slice regular functions}

In this section we summarize the basic definitions from \cite{gs}, \cite{advances},
 that we need in the
sequel.

\pn
Let $\hh$ be the real associative algebra of quaternions
with respect to the basis $\{1, i,j,k \}$
satisfying the relations
$$
i^2=j^2=k^2=-1,
 ij =-ji =k,\
$$
$$
jk =-kj =i ,
 \  ki =-ik =j .
$$
We will denote a quaternion $q$ as $q=x_0+ix_1+jx_2+kx_3$,
$x_i\in \mathbb{R}$, its conjugate as
$\bar q=x_0-ix_1-jx_2-kx_3$, and we will write $|q|^2=q\ov q$.
\noindent

\noindent
Let $\mathbb{S}$ be the sphere of purely imaginary unit quaternions, i.e.
$$
\mathbb{S}=\{ q=ix_1+jx_2+kx_3\ |\  x_1^2+x_2^2+x_3^2=1\}.
$$
\begin{definition} Let $U\subseteq\hh$ be an open set and let
$f:\ U\to\hh$ be a real differentiable function. Let $I\in\mathbb{S}$ and let $f_I$ be
the restriction of $f$ to the complex line $L_I := \mathbb{R}+I\mathbb{R}$
 passing through $1$ and $I$.
We say that $f$ is a slice left regular function if for every $I\in\mathbb{S}$
$$
\frac{1}{2}\left(\frac{\partial }{\partial x}+I\frac{\partial }{\partial y}\right)f_I(x+Iy)=0,
$$
and we say it is slice right regular if for every $I\in\mathbb{S}$
$$
\frac{1}{2}\left(\frac{\partial }{\partial x}+\frac{\partial }{\partial y}I\right)f_I(x+Iy)=0.
$$
\end{definition}

\begin{remark}{\rm Slice left regular functions on $U\subseteq\hh$
form a right vector space $\Er (U)$ and
slice right regular functions on $U\subseteq\hh$ form a left vector space.
 It is not true, in general, that the product of two
regular functions is regular.}
\end{remark}

\noindent
In fact, every slice regular function can be represented as a power series, \cite{advances}:
\begin{theorem} If $B(0,R)$ is a ball centered in the origin with radius $R>0$ and
$f:\ B:\to\hh$ is a slice left regular function, then $f$ has a series expansion of the form
$$
f(q)=\sum_{n=0}^{+\infty} q^n\frac{1}{n!}\frac{\pp^nf}{\pp x^n}(0)
$$
converging on $B$. Analogously, if $f$ is slice right regular it can be expanded as
$$
f(q)=\sum_{n=0}^{+\infty} \frac{1}{n!}\frac{\pp^nf}{\pp x^n}(0) q^n.
$$
\end{theorem}

\begin{remark}{\rm An analogue statement holds for regular functions in an open ball centered in $p_0\in \mathbb{R}$.}
\end{remark}

\begin{definition}
Let $E$ be a bilateral quaternionic Banach space.
A function $f:\ \hh\to E$ is said to be
slice left regular if there exists an open ball $B(0,r)\subset \hh$ and a sequence $\{a_n\}\subset E$ of elements of the Banach space
such that, for every point $q\in B(0,r)$ the function $f(q)$
can be represented by the following  series
\begin{equation}\label{(2.14)}
f(q)=\sum_{n=0}^{+\infty} q^n a_n,
\ \ \ \ \ \ \ \ \    a_n\in E,
\end{equation}
converging in the norm of $E$ for any $q$ such that $|q|< r$. Analogously, $f$ is said to be slice right regular if it can be expanded as
\begin{equation}\label{(2.15)}
f(q)=\sum_{n=0}^{+\infty} a_n q^n ,
\ \ \ \ \ \ \ \ \    a_n\in E,
\end{equation}
\end{definition}

\begin{remark}{\rm From now on we will not specify whether we are considering left or right regular
functions, since the context will clarify it.}
\end{remark}

\begin{definition} Let $V$ be a bilateral vector space on $\mathbb{H}$.
A map $T: V\to V$ is said to be right linear if
$$
 T(u+v)=T(u)+T(v),
$$
$$
T(us)=T(u)s,
$$
for all  $s\in\mathbb{H}$ and  all $u,v\in V$.
\end{definition}
\begin{remark} {\rm
The set End$(V)$ of right linear maps on
$V$ is a ring with respect to the usual sum of operators and with respect to
the composition of operators defined in the usual way:
for any two operators $T,S\in {\rm End}(V)$ we have
$$
(TS)(u)=T(S(u)),\qquad \forall u\in V.
$$
In particular, we have the identity operator
$\mathcal{I}(u)=u$, for all $u\in V$ and
setting $T^0=\mathcal{I}$ we can define powers of a given operator $T\in {\rm End}(V)$:
$T^n=T(T^{n-1})$ for any $n\in \mathbb{N}$. An operator $T$ is said to be invertible
if there exists $S$ such that $TS=ST=\mathcal{I}$ and we will write $S=T^{-1}$.
The set End$(V)$ is a bilateral vector space on $\hh$ with respect to the products
by a scalar defined by
$(aT)(v):=aT(v)$ and
$(Ta)(v):=T(av)$.
\noindent
}
\end{remark}

\subsection{ The case of bounded operators }
\begin{definition} Let $V$ be a bilateral quaternionic Banach space.
We will denote by $\mathcal{B}(V)$ the vector space of all right linear
bounded operators on $V$.
\end{definition}
\noindent
It is easy to verify that
$\mathcal{B}(V)$ is a Banach space endowed with its natural norm.
\begin{definition}  An element $T\in \mathcal{B}(V)$ is said to be invertible if there exists
$T'\in \mathcal{B}(V)$ such that $TT'=T'T=\mathcal{I}$.
\end{definition}
\noindent
It is obvious that the set of all invertible elements of $\mathcal{B}(V)$
 is a group with respect to the composition of operators
defined in $\mathcal{B}(V)$.
\noindent
Let $T$ be a linear quaternionic operator. There are two natural eigenvalue
problems associated to $T$. The first, which one could call the left eigenvalue problem consists
in the solution of equation $T(v)=\lambda v$, and the second, which is called right eigenvalue problem,
and consists in the solution of the equation $T(v)=v\lambda$.

As it was shown in \cite{cgss}, none of them is
 useful to define a functional calculus, as another operator is the one associated
 to the notion of spectrum.
\begin{definition}(The $S$-resolvent operator series)
Let $T\in\mathcal{B}(V)$ and let $s\in \mathbb{H}$. We define the
$S$-resolvent operator series as
\begin{equation}\label{Sresolv}
S^{-1}(s,T):=\sum_{n\geq 0} T^n s^{-1-n}
\end{equation}
for $\|T\|\leq |s|$.
\end{definition}

\begin{theorem}\label{thsmeno} (See \cite{cgss})
Let $T\in\mathcal{B}(V)$ and let $s \in \mathbb{H}$. Assume that
$\overline s$ is such that $T-\overline{s}\mathcal{I}$ is invertible. Then
\begin{equation}\label{STeq}
S(s,T)=(T-\overline s\mathcal{I})^{-1}\, s\, (T-\overline{s}\mathcal{I})-T
\end{equation}
is the inverse of $S^{-1}(s,T)$.
Moreover, we have
\begin{equation}\label{ciao}
\sum_{n\geq 0} T^n s^{-1-n}=-(T^2-2T
Re[s]+|s|^2\mathcal{I})^{-1}(T-\overline{s}\mathcal{I}),
\end{equation}
for $\|T\|\leq |s|$.
\end{theorem}

\begin{definition}(The $S$-resolvent operator)
Let $T\in\mathcal{B}(V)$ and let $s\in \mathbb{H}$.
We define the $S$-resolvent operator as
\begin{equation}\label{Sresolvoperator}
S^{-1}(s,T):=-(T^2-2Re[s] T+|s|^2\mathcal{I})^{-1}(T-\overline{s}\mathcal{I}).
\end{equation}
\end{definition}

\begin{definition}(The $S$-spectrum)
Let $T:V\to V$ be a linear quaternionic operator on the Banach space $V$.
We define the $S$-spectrum $\sigma_S(T)$ of $T$ related to the $S$-resolvent
 operator (\ref{Sresolvoperator}) as:
$$
\sigma_S(T)=\{ s\in \mathbb{H}\ \ :\ \ T^2-2 \ Re [s]T+|s|^2\mathcal{I}\ \ \
{\it is\ not\  invertible}\}.
$$
\end{definition}
\noindent
The following theorem summarizes the properties of the S-spectrum (see \cite{cgss}):
\begin{theorem}\label{struttspettrale}  Let $T\in\mathcal{B}(V)$. Then: \\
1.
The $S$-spectrum $\sigma_S (T)$  is a compact nonempty set contained in
$\{\ s\in \mathbb{H}\ :\  |s|\leq \|T\| \ \}$.
\\
2. If  $p=p_0+p_1I \in \mathbb{H}\setminus \mathbb{R}$ is an $S$-eigenvalue of $T$
 with $I\in \mathbb{S}$ and $p_0, p_1 \in \rr$,
then all the elements of the sphere
$p_0+p_1\mathbb{S}$ are $S$-eigenvalues of $T$.
\par\noindent
Therefore, the  $S$-spectrum $\sigma_S (T)$ is a union  of real points and 2-spheres.
\end{theorem}

Let $S^3=\{q \in \mathbb{H} : |q|=1\}$ denote the unit sphere of $\mathbb{H}$. For any set $A\subseteq \mathbb{H}$, let us define the {\it circularization of} $A$ as  set 
\begin{displaymath}
circ(A):=\bigcup_{x+yI \in A} x+yS^3.
\end{displaymath}
 The following definition will identify an important class of open sets.

\begin{definition} \label{def3.9aa}
Let $T\in\mathcal{B}(V)$.
Let $U \subset \hh$ be an open set such that
\begin{itemize}
\item[(i)] $\pp (U\cap L_I)$ is union of a finite number of
rectifiable Jordan curves  for every $I\in\mathbb{S}$,
\item[(ii) ] $U$ contains the circularization of the $S$-spectrum $\sigma_S(T)$.
\end{itemize}
A function  $f:\ \hh \to \hh$
is said to be locally regular on  $\sigma_S(T)$ if there exists
an open set  $U \subset \hh$, as above,
on which $f$ is regular.
\par\noindent
We will denote
by $\mathcal{R}_{\sigma_S(T)}$ the set of locally regular functions
on $\sigma_S (T)$.

\end{definition}
\begin{theorem}\label{indipdau}
Let $T\in\mathcal{B}(V)$ 
 and  $f\in {\cal R}_{\sigma_S(T)}$.
  Let $U\subset\hh $ be an open set  as in Definition \ref{def3.9aa}  and let $U_I=U\cap L_I$ for $I\in \mathbb{S}$.  Then the integral
\begin{equation}\label{integ311}
{{1}\over{2\pi }} \int_{\partial U_I } S^{-1} (s,T)\  ds_I \ f(s)
\end{equation}
does not depend on the choice of the imaginary unit $I$ and on the
open set $U$.
\end{theorem}

\noindent
The preceding result allows to give the following definition which
offers a new functional calculus:
\begin{definition} \label{fdiT}
Let
 $T\in\mathcal{B}(V)$ and  $f\in {\cal R}_{\sigma_S(T)}$.
Let $U\subset \hh$ be an open set  as in Definition \ref{def3.9aa},
and set $U_I=U\cap L_I$ for $I\in\mathbb{S}$.
We define
\begin{equation}\label{FC}
f(T)= {{1}\over{2\pi }} \int_{\partial U_I } S^{-1} (s,T)\  ds_I \
f(s).
\end{equation}
\end{definition}

\section{Functional calculus for unbounded operators }
We note that if $T$ is a closed operator, the series $\sum_{n\geq 0}T^ns^{-1-n}$ does not
converge.
To overcome this difficulty,
we observe that the right hand side of formula (\ref{ciao}) contains
the continuous operator $(T^2-2T Re[s]+|s|^2\mathcal{I})^{-1}$.
From an heuristical point of view, the composition
$(T^2-2T Re[s]+|s|^2\mathcal{I})^{-1}(T-\overline{s}\mathcal{I})$ gives a bounded operator,
for suitable $s$.
We will consider closed and densely defined operators.

\begin{definition}\label{defschaa}
Let $V$ be a quaternionic Banach space. 
We consider the linear closed densely defined operator
$
T :{\cal D}(T)\subset V\to V
$
where ${\cal D}(T)$ denotes the domain of $T$.
Let us assume that
\begin{itemize}
\item[1)]
${\cal D}(T)$ is dense in $V$,
\item[2)]
$T-\overline{s}\mathcal{I}$  is densely defined in $V$,
\item[3)]
${\cal D}(T^2)\subset {\cal D}(T)$  is dense in $V$,
\item[4)]
$T^2-2T
Re[s]+|s|^2\mathcal{I}$ is one-to-one with range $V$.
\end{itemize}

 The $S$-resolvent operator is defined by
\begin{equation}\label{resochiaa}
S^{-1}(s,T)=-(T^2-2T
Re[s]+|s|^2\mathcal{I})^{-1}(T-\overline{s}\mathcal{I}).
\end{equation}
\end{definition}

\begin{remark}{\rm
We observe that, in principle, it is necessary also the following assumption:
\begin{itemize}
\item[5)] the operator $(T^2-2T Re[s]+|s|^2\mathcal{I})^{-1}(T-\overline{s}\mathcal{I})$ is the restriction to the dense subspace
${\cal D}(T)$ of $V$ of a bounded linear operator.
\end{itemize}
However this assumption is automatically fulfilled since it follows from the identity
$$
(T^2-2T Re[s]+|s|^2\mathcal{I})^{-1}(T-\overline{s}\mathcal{I})
=
T(T^2-2T Re[s]+|s|^2\mathcal{I})^{-1}-
(T^2-2T Re[s]+|s|^2\mathcal{I})^{-1}\overline{s}\mathcal{I},
$$
which is a consequence of
$$
(T^2-2T Re[s]+|s|^2\mathcal{I})^{-1}T=T(T^2-2T Re[s]+|s|^2\mathcal{I})^{-1},
$$
that can be easily verified applying on the left to both sides the operator  $T^2-2T Re[s]+|s|^2\mathcal{I}$.
}
\end{remark}

\begin{definition}\label{Sresolandspecset}
Let $T:{\cal D}(T)\subset V\to V $ be a linear closed densely defined operator as in Definition
\ref{defschaa}.
 We define the $S$-resolvent set of $T$ to be the set
\begin{equation}\label{resosetpoint}
\rho_S(T)=\{ s\in \hh  \ such\ that \ S^{-1}(s,T)
\ exists \ and \ it \ is \  in \ \mathcal{B}(V) \}.
\end{equation}
We define the $S$- spectrum of $T$ as the set
\begin{equation}\label{specsetpoint}
\sigma_S(T)=\hh\setminus \rho_S(T).
\end{equation}
\end{definition}

\begin{theorem}(Structure of the spectrum) \label{strucspec}
 Let $T:V\to V$ be a closed operator
 such that $\sigma_S(T)\not=0$.
 If  $p=p_0+p_1I \in \mathbb{H}\setminus \mathbb{R}$ is an $S$-eigenvalue of $T$
 with $I\in \mathbb{S}$ and $p_0, p_1 \in \rr$,
then all the elements of the sphere
$p_0+p_1\mathbb{S}$ are $S$-eigenvalues of $T$.
The  $S$-spectrum $\sigma_S (T)$ is a union  of real points and 2-spheres.
\end{theorem}
\begin{proof}
The proof is analogous to the one of the bounded case. It immediately  follows from
the structure of the  $S$-eigenvalue equation
 $(T^2-2T Re[s]+|s|^2\mathcal{I})v=0$.
\end{proof}
\begin{theorem}
Let $V$ be a quaternionic Banach space let $T$ be a closed linear quaternionic operator on $V$.
Let $s \in \rho_S(T)$. Then the $S$-resolvent
operator defined in (\ref{Sresolvoperator}) satisfies the equation
$$
S^{-1}(s,T)s-TS^{-1}(s,T)=\mathcal{I}.
$$
\end{theorem}
\begin{proof}
It follows by direct computation. Indeed, replacing (\ref{Sresolvoperator})
in the above equation we have
\begin{equation}\label{serve21}
-(T^2-2Re[s] T+|s|^2\mathcal{I})^{-1}(T-\overline{s}\mathcal{I})s
+T(T^2-2Re[s] T+|s|^2\mathcal{I})^{-1}(T-\overline{s}\mathcal{I})=\mathcal{I}.
\end{equation}
Observe that $T(T^2-2Re[s] T+|s|^2\mathcal{I})^{-1}(T-\overline{s}\mathcal{I})$
 is a bounded operator because it is the sum of two bounded operators.
 Applying now $(T^2-2Re[s] T+|s|^2\mathcal{I})$ to both hands sides of (\ref{serve21}),
we get
$$
-(T-\overline{s}\mathcal{I})s+(T^2-2Re[s] T+|s|^2\mathcal{I})T
(T^2-2Re[s] T+|s|^2\mathcal{I})^{-1}(T-\overline{s}\mathcal{I})=T^2-2Re[s] T+|s|^2\mathcal{I}.
$$
Since $T$ and $T^2-2Re[s] T+|s|^2\mathcal{I}$ commute, we obtain the identity
$$
-(T-\overline{s}\mathcal{I})s+T
(T-\overline{s}\mathcal{I})=T^2-2Re[s] T+|s|^2\mathcal{I}
$$
which proves the statement.
\end{proof}
\begin{definition} Let $s\in \rho_S(T)$. The equation
\begin{equation}\label{resolunbuequa}
S^{-1}(s,T)s-TS^{-1}(s,T)=\mathcal{I}
\end{equation}
will be called the $S$-resolvent equation.
\end{definition}

\begin{theorem}\label{bounSres}
Let $T$ be a closed linear operator on a quaternionic Banach space $V$. Assume
that  $s\in \rho_S(T)$.
Then the $S$-resolvent operator
 can be represented by
\begin{equation}\label{chiusoS}
S^{-1}(s,T)=\sum_{n\geq 0}(Re[s]\mathcal{I}-T)^{-n-1}(Re[s]-s)^n
\end{equation}
if and only if
\begin{equation}\label{Imesse}
|Im[s]|\ \|( Re[s]\mathcal{I}-T)^{-1}\|<1.
\end{equation}
\end{theorem}
\begin{proof}
The equality
$$
(T^2-2T Re[s]+|s|^2\mathcal{I})^{-1}
=
\Big( (T- Re[s]\mathcal{I})^2+|Im[s]|^2\mathcal{I}\Big)^{-1}
$$
$$
=
\Big( (T- Re[s]\mathcal{I})^2 (\mathcal{I}+|Im[s]|^2(T- Re[s]\mathcal{I})^{-2})\Big)^{-1}
$$
$$
=
  \left(\mathcal{I}+|Im[s]|^2(T- Re[s]\mathcal{I})^{-2}\right)^{-1}(T- Re[s]\mathcal{I})^{-2}
$$
$$
=\sum_{n\geq 0}(-1)^n|Im[s]|^{2n}(T- Re[s]\mathcal{I})^{-2n}(T- Re[s]\mathcal{I})^{-2},
$$
yields
$$
S^{-1}(s,T)=
\sum_{n\geq 0}(-1)^{n+1}
|Im[s]|^{2n}(T- Re[s]\mathcal{I})^{-2n-1}
(\mathcal{I}+(T- Re[s]\mathcal{I})^{-1}Im[s])
$$
$$
= \sum_{n\geq 0}(-1)^{n}
|Im[s]|^{2n}(Re[s]\mathcal{I}-T)^{-2n-1}
(\mathcal{I}+(Re[s]\mathcal{I}-T)^{-1}( Re[s] - s))
$$
and since
$$
(-1)^n| Im[s]|^{2n}=( Re[s]-s)^{2n}
$$
we have
$$
S^{-1}(s,T)=
\sum_{n\geq 0}( Re[s]\mathcal{I}-T)^{-2n-1}( Re[s]-s)^{2n}+
\sum_{n\geq 0}( Re[s]\mathcal{I}-T)^{-2n-2}( Re[s]-s)^{2n+1}
$$
$$
=
( Re[s]\mathcal{I}-T)^{-1}
\sum_{n\geq 0}( Re[s]\mathcal{I}-T)^{-n}( Re[s]-s)^{n}
$$
which converges in $\mathcal{B}(V)$ if and only if (\ref{Imesse}) holds.
\end{proof}


Let $V$ be a quaternionic Banach space 
and let $T:{\cal D}(T)\subset V\to V$
be a linear operators.  
If  at least one of its components
is an unbounded operator then its resolvent  is not defined at infinity.
It is therefore natural to consider closed operators $T$ for which the resolvent
$S^{-1}(s,T)$ is not defined at infinity and to define the extended spectrum as
$$
\overline{\sigma}_S(T):=\sigma_S(T)\cup \{\infty\}.
$$
Let us consider $\overline{\hh}=\hh\cup\{\infty\}$ endowed
with the natural topology: a set is open if and only if it is union of open
discs $D(q,r)$ with center at points in $q\in\hh$ and radius $r$, for some $r$, and/or
union of sets the form $\{q\in\hh \ |\ |q|>r\}\cup\{\infty\}=D'(\infty,r)\cup\{\infty\}$, for some $r$.
\begin{definition}
We say that  $f$ is  regular function at $\infty$ if $f(q)$ is
an regular function in a set $D'(\infty,r)$ and
$\lim_{q\to\infty}f(q)$ exists and it is finite. We define $f(\infty)$ to be the
value of this limit.
\end{definition}

\begin{definition}\label{def3.9seconda}
Let $T :{\cal D}(T)\subset V\to V$ be a linear closed operator as in Definition \ref{defschaa} .
Let  $U \subset \hh$ be an open set such that
\begin{itemize}
\item[(i)] $\pp (U\cap L_I)$ is union of a finite number of rectifiable Jordan curves for every $I\in\mathbb{S}$,
\item[(ii)] $U$ contains the circularization of the $S$-spectrum $\sigma_S(T)$.
\end{itemize}
A function  $f$  is said to be locally regular on  $\overline{\sigma}_S(T)$
if it is regular an open set  $U \subset \hh$ as above and at infinity.
\par\noindent
We will denote
by $\mathcal{R}_{\overline{\sigma}_S(T)}$ the set of locally regular functions
on $\overline{\sigma}_S(T)$.
\end{definition}
Consider $\alpha\in\hh$ and the homeomorphism
$$
\Phi :\overline{\hh}\to \overline{\hh}  \ \ \ {\rm for}\ \ \ \alpha\in \hh
$$
 defined by
$$
 p=\Phi(s)=(s-\alpha)^{-1}, \ \ \Phi(\infty)=0,\ \ \ \Phi(\alpha)=\infty.
$$
\begin{definition}
Let $T:{\cal D}(T)\to V$ be a linear closed operator as in Definition \ref{defschaa} with
$\rho_S(T)\cap \mathbb{R}\neq\emptyset$ and suppose that  $f\in \mathcal{R}_{\overline{\sigma}_S(T)}$.
Let us consider
$$
\phi(p):=f(\Phi^{-1}(p))
$$
and the operator
$$
A:=(T-k\mathcal{I})^{-1},\ \ {\it for\ some}\ \  k\in \rho_S(T)\cap \mathbb{R}.
$$
We define
\begin{equation}\label{fdit}
f(T)=\phi(A).
\end{equation}

\end{definition}

\begin{remark}{\rm Observe that, if $\alpha=k\in\rr$, we have that:

i) the function $\phi$ is regular  because it is the
composition
of the function $f$ which is regular and  $\Phi^{-1}(p)=p^{-1}+k$ which is regular
with real coefficients;

ii) in the case $k\in \rho_S(T)\cap \mathbb{R}$  we  have that
$(T-k\mathcal{I})^{-1}=-S^{-1}(k, T)$.
}
\end{remark}
\par\noindent
To prove the fundamental Theorem \ref{fundam} we need the following identities.
\par\noindent
\begin{lemma} Let $s$, $p\in \hh$ and $k\in \mathbb{R}$ such that $p=(s-k)^{-1}$.
Then the following identities hold
\begin{equation}\label{identska}
s_0|p|^2=k|p|^2+p_0,
\end{equation}
\begin{equation}\label{identskb}
 |p|^2|s|^2=k^2|p|^2+2p_0k+1,
\end{equation}
\begin{equation}\label{lambpicc}
(2k-2s_0+\overline{p}^{-1})\frac{\overline{p}}{|p|^2}=-{p^{-2}},
\end{equation}
\begin{equation}\label{lambgrand}
[\overline{s}+(|s|^2-k^2-k\overline{p}^{-1})
(2k-2s_0+\overline{p}^{-1})^{-1}](2k-2s_0+\overline{p}^{-1})\overline{p}=0.
\end{equation}
\end{lemma}
\begin{proof}
Identity (\ref{identska}) follows from 
$$
Re[s-k]=Re[p^{-1}]=Re[\overline{p}|p|^{-2}]
$$
from which we have
$$
s_0-k=p_0|p|^{-2}.
$$
Identity (\ref{identskb}) follows from the chain of identities
$$
|s|^2=s\overline{s}=(k+p^{-1})\overline{(k+p^{-1})}=(k+p^{-1})(k+\overline{p}^{-1})
$$
$$
=
k^2+k(p^{-1}+\overline{p}^{-1})+p^{-1}\overline{p}^{-1}
=k^2+k\frac{2p_0}{|p|^{2}}+\frac{1}{|p|^{2}}.
$$
To prove (\ref{lambpicc}) we consider the chain of identities
$$
(2k-2s_0+\overline{p}^{-1})\frac{\overline{p}}{|p|^2}
=
(2k-2s_0+\overline{p}^{-1}){p}^{-1}
$$
$$
=
(2k-2s_0+\overline{s}-k)(s-k)=-(s-k)^2=-{p^{-2}}.
$$
Finally (\ref{lambgrand}) follows from
$$
[\overline{s}+(|s|^2-k^2-k\overline{p}^{-1})
(2k-2s_0+\overline{p}^{-1})^{-1}](2k-2s_0+\overline{p}^{-1})\overline{p}
$$
$$
=
\overline{s}(2(k-s_0)+\overline{p}^{-1})\overline{p}
+(|s|^2-k^2-k\overline{p}^{-1})\overline{p}
$$
now using (\ref{identska}) and (\ref{identskb}) we get
$$
\overline{s}(2(k-s_0)+\overline{p}^{-1})\overline{p}
+(|s|^2-k^2-k\overline{p}^{-1})\overline{p}
$$
$$
=\overline{s}\Big(-2\frac{p_0}{|p|^2}\overline{p}+1\Big)
+\frac{2p_0k+1}{|p|^2}\overline{p}-k
$$
$$
=(k+(\overline{p})^{-1})\Big(-2\frac{p_0}{|p|^2}\overline{p}+1\Big)
+\frac{2p_0k+1}{|p|^2}\overline{p}-k
$$
$$
=-2p_0\overline{p}^{-1}p^{-1}+\overline{p}^{-1}+p^{-1}=0.
$$
\end{proof}

\begin{theorem}\label{fundam}
If $k\in \rho_S(T)\cap\rr\not=\emptyset$ and $\Phi$, $\phi$ are as above, then $\Phi(\overline{\sigma}_S(T))=\sigma_S(A)$ and the relation
$\phi(p):=f(\Phi^{-1}(p))$ determines a one-to-one correspondence  between $f\in \mathcal{R}_{\overline{\sigma}_S(T)}$ and $\phi\in \mathcal{R}_{\overline{\sigma}_S(A)}$.
\end{theorem}
\begin{proof} First we consider the case
 $p\in \sigma_S(A)$ and $p\not= 0$. Recall that
$$
S^{-1}(p,A)=-(A^2-2A Re[p]+|p|^2\mathcal{I})^{-1}(A-\overline{p}\mathcal{I}),
$$
from which we obtain
$$
(A^2-2A Re[p]+|p|^2\mathcal{I})S^{-1}(p,A)=-(A-\overline{p}\mathcal{I}).
$$
Let us apply the operator $A^{-2}$ on the left to get
$$
(\mathcal{I}-2 Re[p] A^{-1}+A^{-2}|p|^2)S^{-1}(p,A)=-(A^{-1}-A^{-2}\overline{p}).
$$
Now we use the relations
\begin{equation}\label{AeAduemeno}
A^{-1}=T-k\mathcal{I},\ \ \ A^{-2}=T^2-2kT+k^2\mathcal{I}
\end{equation}
to get
$$
(\mathcal{I}-2 Re[p] (T-k\mathcal{I})+(T^2-2kT+k^2\mathcal{I})|p|^2)S^{-1}(p,A)
$$
$$
=-(T-k\mathcal{I}-(T^2-2kT+k^2\mathcal{I})\overline{p}).
$$
Using the identities (\ref{identska}) and (\ref{identskb})
we have
$$
(T^2 |p|^2 -2s_0|p|^2T + |s|^2 |p|^2\mathcal{I})
S^{-1}(p,A)=-(T-k\mathcal{I}-(T^2-2kT+k^2\mathcal{I})\overline{p}).
$$
So we get the equalities
$$
S^{-1}(p,A)=-\frac{1}{|p|^2}(T^2  -2s_0T + |s|^2 \mathcal{I})^{-1}
(T-k\mathcal{I}-(T^2-2kT+k^2\mathcal{I})\overline{p})
$$
$$
=-\frac{1}{|p|^2}(T^2  -2s_0T + |s|^2 \mathcal{I})^{-1}
(T{\overline{p}}^{-1}-k{\overline{p}}^{-1}\mathcal{I}-T^2+2kT-k^2\mathcal{I})\overline{p}
$$
$$
=-\frac{1}{|p|^2}(T^2  -2s_0T + |s|^2 \mathcal{I})^{-1}
\Big(-(T^2-2s_0T+|s|^2\mathcal{I})\overline{p}
$$
$$
+( T(2k-2s_0+\overline{p}^{-1})
+(|s|^2-k^2-k\overline{p}^{-1})\mathcal{I})\overline{p}\Big)
$$
$$
=
\frac{\overline{p}}{|p|^2}\mathcal{I}
-\frac{1}{|p|^2}(T^2  -2s_0T + |s|^2 \mathcal{I})^{-1}
( T
+(|s|^2-k^2-k\overline{p}^{-1})(2k-2s_0
+\overline{p}^{-1})^{-1}\mathcal{I})(2k-2s_0+\overline{p}^{-1})\overline{p}.
$$
With some calculation we get
$$
S^{-1}(p,A)=p^{-1}\mathcal{I}
$$
$$
-\frac{1}{|p|^2}(T^2  -2s_0T + |s|^2 \mathcal{I})^{-1}
( T -\overline{s}\mathcal{I}+
[\overline{s}+(|s|^2-k^2-k\overline{p}^{-1})
(2k-2s_0+\overline{p}^{-1})^{-1}] \mathcal{I})(2k-2s_0+\overline{p}^{-1})\overline{p}
$$
and also
$$
S^{-1}(p,A)=
p^{-1}\mathcal{I}
+S^{-1}(s,T)\lambda
-\frac{1}{|p|^2}(T^2  -2s_0T + |s|^2 \mathcal{I})^{-1}
\Lambda
$$
where we have set
$$
\lambda:=(2k-2s_0+\overline{p}^{-1})\frac{\overline{p}}{|p|^2},
$$
$$
\Lambda:=[\overline{s}+(|s|^2-k^2-k\overline{p}^{-1})
(2k-2s_0+\overline{p}^{-1})^{-1}](2k-2s_0+\overline{p}^{-1})\overline{p}.
$$
Using the identities (\ref{lambpicc}) and (\ref{lambgrand}) we finally get 
$$
S^{-1}(p,A)=
p^{-1}\mathcal{I}
-S^{-1}(s,T){p^{-2}},
$$
but also
\begin{equation}\label{importante}
S^{-1}(s,T)=
p\mathcal{I}
-S^{-1}(p,A)p^2.
\end{equation}
So $p\in \rho_S(A)$, $p\not=0$ then $s\in \rho_S(T)$.
\par\noindent
Now take $s\in \rho_S(T)$ and observe that, with form the definitions of $S^{-1}(s,T)$ and 
of $A$, 
replacing (\ref{AeAduemeno}) in the $S$ resolvent operator $S^{-1}(s,T)$ 
we obtain
$$
S^{-1}(s,T)=-(A^{-2}-2p_0|p|^{-2}A^{-1} +|p|^{-2}\mathcal{I})^{-1}(A^{-1}+k-\overline{s}\mathcal{I}),
$$
where we also have used identity (\ref{identskb}).
By some calculation we finally get: 
$$
S^{-1}(s,T)=-AS^{-1}(p,A)p,
$$
so if  $s\in \rho_S(T)$ then $p\in \rho_S(A)$, $p\not=0$.
\par\noindent
The point $p=0$ belongs to $\sigma_S(A)$
 since $S^{-1}(0,A)=A^{-1}=T-k\mathcal{I}$ is unbounded.
 The last part of the statement is evident from the definition of $\Phi$.
\end{proof}

\begin{theorem}\label{calcfunzformu}
Let $T:{\cal D}(T)\subset V\to V$ be a linear closed operator as in Definition \ref{defschaa} with
$\rho_S(T)\cap \mathbb{R}\neq\emptyset$ and suppose that  $f\in \mathcal{R}_{\overline{\sigma}_S(T)}$.
Then operator $f(T)$ defined in (\ref{fdit}) is independent of $k\in \rho_S(T)\cap \mathbb{R}$.

Let $W$, be an open set such that $\overline{\sigma}_S(T) \subset W$
and let $f$ be an regular function on $W\cup\partial W$.
Set $W_I= W\cap L_I$ for $I\in \mathbb{S}$ be such that
its boundary $\partial W_I$ is positively oriented and consists of a finite number
of rectifiable Jordan curves.
Then
\begin{equation}\label{fditfor}
f(T)=f(\infty)\mathcal{I}+\frac{1}{2\pi} \int_{\partial W_I} S^{-1}(s,T)ds_I f(s).
\end{equation}
\end{theorem}
\begin{proof}
The first part of the statement follows from the validity of formula (\ref{fditfor})
since the integral is independent of $k$.
\par\noindent
Given $k\in\rho_S(T)\cap\mathbb{R}$ and the set $W$ we can assume that
$k\not\in W_I\cup\partial W_I$, $\forall I\in\mathbb{S}$
since otherwise, by the Cauchy theorem,  we can replace
$W$ by $W'$, on which $f$ is regular, such that $k\not\in W_I'\cup\partial W_I'$,
without altering the value of the integral (\ref{fditfor}). Moreover, the integral
(\ref{fditfor}) is independent of the choice of $I\in\mathbb{S}$, thanks to the
structure of the spectrum (see Theorem \ref{strucspec}) and an argument similar to the
one used to prove Theorem \ref{indipdau}.
\par\noindent
 We have that ${\cal V}_I:=\Phi^{-1}(W_I)$ is an open set that contains
$\sigma_S(T)$ and its boundary $\partial {\cal V}_I=\Phi^{-1}(\partial W_I)$
is positively oriented and consists of a finite number
of rectifiable Jordan curves.  Using the relation (\ref{importante}) we have 
$$
  \frac{1}{2\pi} \int_{\partial W_I} S^{-1}(s,T)ds_I f(s)
$$
$$
  = -\frac{1}{2\pi} \int_{\partial {\cal V}_I} \Big(p\mathcal{I}-S^{-1}(p,A)p^2 \Big) p^{-2}dp_I \phi(p)
$$
$$= -\frac{1}{2\pi} \int_{\partial {\cal V}_I} p^{-1}dp_I \phi(p)+\frac{1}{2\pi} \int_{\partial {\cal V}_I} S^{-1}(p,A)dp_I \phi(p)
$$
$$
= -\mathcal{I} \phi(0)+\phi(A)
$$
now by definition $\phi(A)=f(T)$ and $\phi(0)=f(\infty)$  we obtain 
$$
\frac{1}{2\pi} \int_{\partial W_I} S^{-1}(s,T)ds_I f(s)=-\mathcal{I} f(\infty)+f(T).
$$

\end{proof}

\begin{theorem}\label{Tmenouno}
Let $T$ be a closed quaternionic operator with a bounded inverse, and such that
$\rho_S(T)\cap\mathbb{R}\not=\emptyset$. Suppose that (\ref{Imesse}) holds.
Let $V$ be an open set in $L_I$ that contains $L_I\cap \sigma_S(T)$
such that its boundary $\partial V$ consists of the segment $I_R$ of length $2R$ on the
imaginary axis $I$ symmetric with respect to the origin
 and of the semicircle $\gamma_R$ with diameter $I_R$ and $Re[s]>0$.
Let $f$ be slice regular  on $\overline{V}\cup\infty$ such that $f(\infty)=0$ and suppose that
$\int_{ \gamma_R } S^{-1}(s,T)ds_I f(s)\to 0$ in $\mathcal{B}(V)$ as $R\to +\infty$.
Then
\begin{equation}
f(T)=\sum_{n\geq 0}
T^{-n-1}\mathcal{F}_n(f),
\end{equation}
where
\begin{equation}
\mathcal{F}_n(f)=-\frac{1}{2\pi} \int_{I\mathbb{R}} (Im[s])^{n} ds_I f(s),
\end{equation}
when the integral converges.
\end{theorem}
\begin{proof}
From the previous theorem, we immediately obtain
$$
f(T)=\frac{1}{2\pi} \int_{I_R\cup \gamma_R } S^{-1}(s,T)ds_I f(s)
$$
$$
=\frac{1}{2\pi} \int_{I_R } S^{-1}(s,T)ds_I f(s)
+\frac{1}{2\pi} \int_{ \gamma_R } S^{-1}(s,T)ds_I f(s)
:=A_1(R, f)+A_2(R, f)
$$
By the hypotheses on $f$ we have that
$
\|A_2(R, f)\|\to 0
$
as $R\to \infty$, and therefore
$
f(T)=A_1(\infty, f).
$
In this case  the $S$-resolvent representation
(\ref{chiusoS}) implies
$$
S^{-1}(s,T)=\sum_{n\geq 0}(-1)^{n+1}T^{-n-1}(-1)^n(Im[s])^n=-\sum_{n\geq 0}T^{-n-1}(Im[s])^n.
$$
The statement follows from the definition of $A_1(\infty, f)$ and $\mathcal{F}_n(f)$.
\end{proof}
\begin{remark}{\rm  The Cauchy theorem shows that any set $V$ whose boundary consists of a finite number
of regular curves can be assumed to satisfy the conditions in the statement.}
\end{remark}

\end{document}